# Рандомизация и разреженность в задачах huge-scale оптимизации на примере работы метода зеркального спуска

*Аникин А.С. (ИДСТУ СО РАН)*
*Гасников А.В. (ИППИ РАН; ПреМоЛаб МФТИ)*
*Горнов А.Ю. (ИДСТУ СО РАН)*

**Аннотация**
В работе исследуются различные рандомизации метода зеркального спуска для задач huge-scale оптимизации с разреженной структурой. В качестве одного из примеров приложения приводится задача PageRank.
**Ключевые слова:** huge-scale оптимизация, рандомизация, метод зеркального спуска, разреженность, оценки вероятностей больших уклонений, PageRank.

## 1. Введение

В недавнем цикле работ Ю.Е. Нестерова с соавторами [1 – 4] был введен класс задач huge-scale оптимизации (задачи выпуклой оптимизации, для которых размерность прямого и(или) двойственного пространства не меньше десятков миллионов), и исследовалась роль разреженности в таких задачах. Настоящая работа посвящена изучению конкретного (но, пожалуй, наиболее важного) метода решения таких задач – метода зеркального спуска (МЗС). Это метод представляет собой обобщение метода проекции градиента. Он был предложен в конце 70-х годов прошлого века А.С. Немировским [5]. С тех пор метод получил повсеместное распространение для решения задач больших размерностей, прежде всего, в связи с идей рандомизации. Метод оказался слабо чувствительным к замене настоящего градиента его несмещенной оценкой. Это обстоятельство активно используется на практике, поскольку построить несмещенную оценку в ряде случаев удается намного дешевле, чем посчитать градиент. Как правило, выгода получается пропорциональной размерности пространства, в котором происходит оптимизация. Однако для разреженных задач тут возникают нюансы, связанные с проработкой правильного сочетания рандомизации и разреженности. Настоящая работа посвящена изучению такого сочетания для МЗС. Данная работа является продолжением работы [3], в которой уже был рассмотрен один сюжет на эту тему (метод Григориадиса–Хачияна для задачи PageRank).

Структура статьи следующая. В п. 2 мы описываем рандомизированый МЗС. Стоит обратить внимание на случай, когда оптимизация происходит на неограниченном множестве, и при этом выводятся оценки вероятностей больших уклонений. Тут есть некоторые тонкости, проработка которых делает п. 2 не просто вводным материалом для последующего изложения, но и представляющим самостоятельный интерес. В п. 3 результаты п. 2 переносятся на задачи с функциональными ограничениями, на которые мы не умеем эффективно проектироваться. В детерминированном случае (когда используется обычный градиент) такого типа задачи рассматривались достаточно давно. Разработаны эффективные способы их редуцирования к задаче с простыми ограничениями. Предлагались различные эффективные методы (Поляк–Шор–Немировский–Нестеров). Однако в случае, когда вместо градиента мы используем его несмещенную оценку (для функционала и ограничений) нам не известны оценки, поэтому в п. 3 приводится соответствующее обобщение МЗС и устанавливаются необходимые в дальнейшем оценки. В п. 4 на основе теоретических заготовок пп. 2, 3 мы описываем класс разреженных задач (обобщающих задачу PageRank), для которых удается за счет рандомизации получить дополнительную выгоду.

## 2. Рандомизированный метод зеркального спуска

Рассмотрим задачу выпуклой оптимизации
$$f(x) \to \min_{x \in Q}. \qquad (1)$$
Под решением этой задачи будем понимать такой $\bar{x}^N \in Q \subseteq \mathbb{R}^n$, что с вероятностью $\geq 1 - \sigma$ имеет место неравенство
$$f(\bar{x}^N) - f_* \leq \varepsilon,$$
где $f_* = f(x_*)$ – оптимальное значение функционала в задаче (1), $x_* \in Q$ – решение задачи (1). На каждой итерации $k = 1, ..., N$ нам доступен стохастический (суб-)градиент $\nabla_x f_k(x, \xi^k)$ в одной, выбранной нами (методом), точке $x^k$.

Опишем метод зеркального спуска (МЗС) для решения задачи (1) (мы в основном будем следовать работам [5, 6]). Введем норму $\|\ \|$ в прямом пространстве (сопряженную норму будем обозначать $\|\ \|_*$) и прокс-функцию $d(x)$ сильно выпуклую относительно этой нормы, с константой сильной выпуклости $\geq 1$. Выберем точку старта
$$x^1 = \arg\min_{x \in Q} d(x),$$
считаем, что
$$d(x^1) = 0, \ \nabla d(x^1) = 0.$$
Введем брэгмановское "расстояние"
$$V_x(y) = d(y) - d(x) - \langle \nabla d(x), y - x \rangle.$$
Определим "размер" решения
$$d(x_*) = V_{x^1}(x_*) = R^2.$$
Определим оператор "проектирования" согласно этому расстоянию
$$\mathrm{Mirr}_{x^k}(\mathrm{v}) = \arg\min_{y \in Q} \{\langle \mathrm{v}, y - x^k \rangle + V_{x^k}(y)\}.$$
МЗС для задачи (1) будет иметь вид, см., например, [6]
$$x^{k+1} = \mathrm{Mirr}_{x^k}\left(\alpha \nabla_x f_k(x^k, \xi^k)\right), \ k = 1, ..., N.$$

Будем считать, что $\{\xi^k\}_{k=1}^N$ представляет собой последовательность независимых случайных величин, и для всех $x \in Q$ имеют место условия ($k = 1, ..., N$)
1. $E_{\xi^k}\left[\nabla_x f_k(x, \xi^k)\right] = \nabla f(x);$
2. $E_{\xi^k}\left[\|\nabla_x f_k(x, \xi^k)\|_*^2\right] \leq M^2.$

В ряде случаев нам также понадобится более сильное условие

3. $\|\nabla_x f_k(x, \xi^k)\|_*^2 \leq \tilde{M}^2$ почти наверное по $\xi^k$.

При выполнении условия 1 для любого $u \in Q$, $k = 1, ..., N$ имеет место неравенство, см., например, [6]
$$\alpha \langle \nabla_x f_k(x^k, \xi^k), x^k - u \rangle \leq \frac{\alpha^2}{2} \|\nabla_x f_k(x^k, \xi^k)\|_*^2 + V_{x^k}(u) - V_{x^{k+1}}(u).$$
Это неравенство несложно получить в случае евклидовой прокс-структуры [7]

$$d(x) = \|x\|_2^2 / 2, \quad V_x(y) = \|y - x\|_2^2 / 2.$$

В этом случае МЗС для задачи (1) есть просто вариант обычного метода проекции градиента (см. примеры 1, 2 ниже).

Разделим сначала выписанное неравенство на $\alpha$ и возьмем условное математическое ожидание $E_{\xi^k}\left[\,\cdot\,\big|\Xi^{k-1}\right]$ ($\Xi^{k-1}$ — сигма алгебра, порожденная $\xi^1$, …, $\xi^{k-1}$), затем просуммируем то, что получится по $k = 1, ..., N$, используя условие 1. Затем возьмем от того, что получилось при суммировании, полное математическое ожидание, учитывая условие 2. В итоге, выбрав $u = x_*$, и определив

$$\bar{x}^N = \frac{1}{N} \sum_{k=1}^{N} x^k,$$

получим

$$N \cdot \left(E\left[f(\bar{x}^N)\right] - f_*\right) \le \frac{V_{x^1}(x_*)}{\alpha} - \frac{E\left[V_{x^{N+1}}(x_*)\right]}{\alpha} + \frac{1}{2} M^2 \alpha N \le$$
$$\le \frac{R^2}{\alpha} + \frac{1}{2} M^2 \alpha N.$$

Выбирая[1]

$$\alpha = \frac{R}{M}\sqrt{\frac{2}{N}},$$

получим

$$E\left[f(\bar{x}^N)\right] - f_* \le MR\sqrt{\frac{2}{N}}. \qquad (2)$$

Заметим, что в детерминированном случае вместо $\bar{x}^N$ можно брать
$$\breve{x}^N = \arg\min_{k=1,...,N} f(x^k).$$

Немного более аккуратные рассуждения (использующие неравенство Азума–Хефдинга) с

$$\alpha = \frac{R}{\tilde{M}}\sqrt{\frac{2}{N}}$$

позволяют уточнить оценку (2) следующим образом (см., например, [8, 9]):

$$f(\bar{x}^N) - f_* \le \tilde{M}\sqrt{\frac{2}{N}}\left(R + 2\tilde{R}\sqrt{\ln(2/\sigma)}\right) \qquad (3)$$

с вероятностью $\ge 1 - \sigma$, где

$$\tilde{R} = \sup_{x \in \tilde{Q}} \|x - x_*\|,$$
$$\tilde{Q} = \left\{x \in Q : \|x - x_*\|^2 \le 65 R^2 \ln(4N/\sigma)\right\}.$$

Собственно, для справедливости оценки (3) достаточно требовать выполнение условий 1 – 3 лишь на множестве $\tilde{Q} \subseteq Q$. Это замечание существенно, когда рассматриваются неограниченные множества $Q$ (см., например, п. 4).

Оценки (2), (3) являются неулучшаемыми с точностью до мультипликативного числового множителя по $N$ и $\sigma$. Наряду с этим можно обеспечить и их неулучшаемость по размерности пространства $n$, путем "правильного" выбора прокс-функции [5] (такой выбор всегда возможен, и известен для многих важных в приложениях случаев выбора множества $Q$). Собственно прокс-структура (новая степень свободы по сравнению с классическим методом проекции градиента) и вводилась для того, чтобы была возможность обеспечить последнее свойство.

---
[1] Можно получить и адаптивный вариант приводимой далее оценки, для этого потребуется использовать метод двойственных усреднений [6, 7]. Впрочем, в [5] имеется адаптивный вариант МЗС.

В виду того, что мы используем рандомизированный метод и всегда
$$f(\bar{x}^N) - f_* \geq 0,$$
то используя идею амплификации (широко распространенную в Computer Science) можно немного "улучшить" оценку (3). Для этого сначала перепишем оценку (2) в виде
$$E[f(\bar{x}^N)] - f_* \leq \varepsilon,$$
где (здесь $N$, конечно, должно быть натуральным числом, поэтому эту формулу и последующие формулы такого типа надо понимать с точностью до округления к наименьшему натуральному числу, большему написанного)
$$N = \frac{2M^2 R^2}{\varepsilon^2}. \qquad (4)$$

Отсюда по неравенству Маркова
$$P(f(\bar{x}^N) - f_* \geq 2\varepsilon) \leq \frac{E[f(\bar{x}^N)] - f_*}{2\varepsilon} \leq \frac{1}{2}.$$

Можно параллельно (независимо) запустить $\log_2(\sigma^{-1})$ траекторий метода. Обозначим $\bar{x}_{\min}^N$ тот из $\bar{x}^N$ на этих траекториях, который доставляет минимальное значение $f(\bar{x}^N)$. Из выписанного неравенства Маркова получаем, что имеет место неравенство
$$P(f(\bar{x}_{\min}^N) - f_* \geq 2\varepsilon) \leq \sigma.$$

Таким образом, можно не более чем за
$$N = \frac{8M^2 R^2}{\varepsilon^2} \log_2(\sigma^{-1})$$
обращений за стохастическим градиентом и не более чем за $\log_2(\sigma^{-1})$ обращений за значением функции найти решение задачи (1) $\bar{x}^N$ с требуемой точностью $\varepsilon$ и доверительным уровнем $\sigma$.

Как уже отмечалось, во многих приложениях множество $Q$ неограниченно (см., например, п. 4). Поскольку $x_*$ априорно не известно, то это создает проблемы для определения $R$, которое входит в формулу для расчета шага метода
$$\alpha = \frac{R}{M}\sqrt{\frac{2}{N}}.$$

Однако если мы заранее выбрали желаемую точность $\varepsilon$, то с помощью формулы (4) можно выразить шаг следующим образом
$$\alpha = \frac{\varepsilon}{M^2}.$$

Рассмотрим три конкретных примера множества $Q$, которые нам понадобятся в дальнейшем (см. п. 4). В примерах 1, 2 мы не приводим оценки скорости сходимости, поскольку они будут иметь вид (2), (3), т.е. никакой уточняющей информации тут не появляется, в отличие от примера 3.

**Пример 1 (все пространство).** Предположим, что $Q = \mathbb{R}^n$. Выберем
$$\|\ \| = \|\ \|_2,\ d(x) = \frac{1}{2}\|x\|_2^2.$$
Тогда МЗС примет следующий вид ($\alpha = \varepsilon/M^2$, $x^1 = 0$):

$$x^{k+1} = x^k - \alpha \nabla_x f_k\left(x^k, \xi^k\right), \ k = 1,...,N. \ \square$$

**Пример 2 (неотрицательный ортант).** Предположим, что
$$Q = \mathbb{R}^n_+ = \left\{x \in \mathbb{R}^n: \ x \geq 0\right\}.$$

Выберем
$$\| \ \| = \| \ \|_2, \ d(x) = \frac{1}{2}\|x - \bar{x}\|_2^2, \ \bar{x} \in \text{int } Q.$$

Тогда МЗС примет следующий вид ($\alpha = \varepsilon/M^2$, $x^1 = \bar{x}$):
$$x^{k+1} = \left[x^k - \alpha \nabla_x f_k\left(x^k, \xi^k\right)\right]_+ = \max\left\{x^k - \alpha \nabla_x f_k\left(x^k, \xi^k\right), 0\right\}, \ k = 1,...,N,$$

где $\max\{ \ \}$ берется покомпонентно. $\square$

**Пример 3 (симплекс).** Предположим, что
$$Q = S_n(1) = \left\{x \geq 0: \ \sum_{i=1}^n x_i = 1\right\}.$$

Выберем
$$\| \ \| = \| \ \|_1, \ d(x) = \ln n + \sum_{i=1}^n x_i \ln x_i.$$

Тогда МЗС примет следующий вид:
$$x_i^1 = 1/n, \ i = 1,...,n,$$

при $k = 1,...,N$, $i = 1,...,n$
$$x_i^{k+1} = \frac{\exp\left(-\sum_{r=1}^k \alpha \frac{\partial f_r(x^r, \xi^r)}{\partial x_i}\right)}{\sum_{l=1}^n \exp\left(-\sum_{r=1}^k \alpha \frac{\partial f_r(x^r, \xi^r)}{\partial x_l}\right)} = \frac{x_i^k \exp\left(-\alpha \frac{\partial f_k(x^k, \xi^k)}{\partial x_i}\right)}{\sum_{l=1}^n x_l^k \exp\left(-\alpha \frac{\partial f_k(x^k, \xi^k)}{\partial x_l}\right)}.$$

Оценки скорости сходимости будут иметь вид:
$$E\left[f(\bar{x}^N)\right] - f_* \leq M\sqrt{\frac{2\ln n}{N}} \text{ (при } \alpha = M^{-1}\sqrt{2\ln n/N}\text{)};$$

$$f(\bar{x}^N) - f_* \leq \tilde{M}\sqrt{\frac{2}{N}}\left(\sqrt{\ln n} + 4\sqrt{\ln(\sigma^{-1})}\right) \text{ (при } \alpha = \tilde{M}^{-1}\sqrt{2\ln n/N}\text{)}$$

с вероятностью $\geq 1 - \sigma$. $\square$

Представим себе, что задача (1) видоизменилась следующим образом
$$f(x) \to \min_{\substack{g(x) \leq 0 \\ x \in Q}}.$$

Можно ли к ней применить изложенный выше подход, если "проектироваться" на множество (в отличие от $Q$)
$$\{x \in Q: \ g(x) \leq 0\}$$

мы эффективно не умеем? В случае, когда мы знаем оптимальное значение $f_*$, то мы можем свести новую задачу к задаче (1)
$$\min\{f(x) - f_*, g(x)\} \to \min_{x \in Q}.$$

Несложно записать рандомизированный МЗС для такой задачи и применить к ней все, что изложено выше. Однако мы не будем здесь этого делать, поскольку в следующем пункте мы приведем более общий вариант МЗС, который не предполагает, что известно $f_*$.

### 3. Рандомизированный метод зеркального спуска с функциональными ограничениями

Рассмотрим задачу

$$f(x) \to \min_{\substack{g(x) \leq 0 \\ x \in Q}}. \tag{5}$$

Под решением этой задачи будем понимать такой $\bar{x}^N \in Q$, что имеют место неравенства

$$E\left[f\left(\bar{x}^N\right)\right] - f_* \leq \varepsilon_f = \frac{M_f}{M_g}\varepsilon_g, \quad g\left(\bar{x}^N\right) \leq \varepsilon_g,$$

где $f_* = f(x_*)$ – оптимальное значение функционала в задаче (5), $x_* \in Q$ – решение задачи (5). Будем считать, что имеется такая последовательность независимых случайных величин $\{\xi^k\}$ и последовательности $\{\nabla_x f_k(x, \xi^k)\}$, $\{\nabla_x g_k(x, \xi^k)\}$, $k = 1, \ldots, N$, что для всех $x \in Q$ имеют место следующие соотношения (можно считать их выполненными при $x \in \tilde{Q}$, см. п. 2 и последующий текст в этом пункте)

$$E_{\xi^k}\left[\nabla_x f_k(x, \xi^k)\right] = \nabla f(x), \quad E_{\xi^k}\left[\nabla_x g_k(x, \xi^k)\right] = \nabla g(x);$$

$$E_{\xi^k}\left[\left\|\nabla_x f_k(x, \xi^k)\right\|_*^2\right] \leq M_f^2, \quad E_{\xi^k}\left[\left\|\nabla_x g_k(x, \xi^k)\right\|_*^2\right] \leq M_g^2.$$

МЗС для задачи (5) будет иметь вид (см., например, в детерминированном случае двойственный градиентный метод из [10])

$$x^{k+1} = \text{Mirr}_{x^k}\left(h_f \nabla_x f_k(x^k, \xi^k)\right), \text{ если } g(x^k) \leq \varepsilon_g,$$

$$x^{k+1} = \text{Mirr}_{x^k}\left(h_g \nabla_x g_k(x^k, \xi^k)\right), \text{ если } g(x^k) > \varepsilon_g,$$

где $h_g = \varepsilon_g/M_g^2$, $h_f = \varepsilon_g/(M_f M_g)$, $k = 1, \ldots, N$. Обозначим через $I$ множество индексов $k$, для которых $g(x^k) \leq \varepsilon$. Введем также обозначения

$$[N] = \{1, \ldots, N\}, \quad J = [N] \setminus I, \quad N_I = |I|, \quad N_J = |J|, \quad \bar{x}^N = \frac{1}{N_I}\sum_{k \in I} x^k.$$

Тогда имеет место неравенство

$$h_f N_I \cdot \left(E\left[f\left(\bar{x}^N\right)\right] - f_*\right) \leq$$

$$\leq h_f E\left[\sum_{k \in I}\left\langle E_{\xi^k}\left[\nabla_x f_k(x^k, \xi^k)\right], x^k - x_*\right\rangle\right] \leq \frac{h_f^2}{2}\sum_{k \in I} E\left[\left\|\nabla_x f_k(x^k, \xi^k)\right\|_*^2\right] +$$

$$- h_g E\Big[\sum_{k \in J}\underbrace{\left\langle E_{\xi^k}\left[\nabla_x g_k(x^k, \xi^k)\right], x^k - x_*\right\rangle}_{\geq g(x^k) - g(x_*) > \varepsilon_g}\Big] + \frac{h_g^2}{2}\sum_{k \in J} E\left[\left\|\nabla_x g_k(x^k, \xi^k)\right\|_*^2\right] +$$

$$+ \sum_{k \in [N]}\left(E\left[V_{x^k}(x_*)\right] - E\left[V_{x^{k+1}}(x_*)\right]\right) \leq$$

$$\leq \frac{1}{2}h_f^2 M_f^2 N_I - \frac{1}{2M_g^2}\varepsilon_g^2 N_J + E\left[V_{x^1}(x_*)\right] - E\left[V_{x^{N+1}}(x_*)\right] =$$

$$= \frac{1}{2}\left(h_f^2 M_f^2 + \frac{\varepsilon_g^2}{M_g^2}\right)N_I - \frac{1}{2M_g^2}\varepsilon_g^2 N + R^2 - E\left[V_{x^{N+1}}(x_*)\right].$$

Будем считать, что (следует сравнить с формулой (4))

$$N = N(\varepsilon_g) = \frac{2M_g^2 R^2}{\varepsilon_g^2} + 1. \qquad (6)$$

Тогда $E[N_I] \geq 1$ и

$$E\left[f(\bar{x}^N)\right] - f_* \leq \frac{1}{2}\left(h_f M_f^2 + \frac{\varepsilon_g^2}{M_g^2 h_f}\right) = \frac{M_f}{M_g}\varepsilon_g = \varepsilon_f.$$

Соотношение

$$g(\bar{x}^N) \leq \varepsilon_g$$

следует из того, что по построению $g(x^k) \leq \varepsilon_g$, $k \in I$ и из выпуклости функции $g(x)$.

Заметим, что в детерминированном случае вместо $\bar{x}^N$ можно брать

$$\breve{x}^N = \arg\min_{k \in I} f(x^k).$$

Если известно, что для всех $x \in \tilde{Q}$ и почти наверное по $\xi^k$

$$\left\|\nabla_x f_k(x,\xi^k)\right\|_*^2 \leq \tilde{M}^2, \quad \left\|\nabla_x g_k(x,\xi^k)\right\|_*^2 \leq \tilde{M}^2, \quad k = 1,\ldots,N,$$

то для описанного в этом пункте метода (с $\varepsilon = \varepsilon_f = \varepsilon_g$ и $h_f = h_g = \varepsilon/\tilde{M}^2$) вид оценки вероятностей больших уклонений (3) из п. 2 сохранится (оценка получается чуть лучше, чем нижняя оценка из работы [11], когда ограничений-неравенств больше одного, поскольку мы имеем доступ к точному значению $g(x)$)

$$f(\bar{x}^N) - f_* \leq 2\tilde{M}\sqrt{\frac{2}{N}}\left(R + 2\tilde{R}\sqrt{\ln(2/\sigma)}\right),$$

где $N = 2N(\varepsilon)$ (см. формулу (6)). К сожалению, трюк с амплификацией (см. п. 2) здесь уже не проходит в том же виде, как и раньше, поскольку теперь уже нельзя гарантировать

$$f(\bar{x}^N) - f_* \geq 0.$$

Однако если ввести обозначение

$$\varepsilon_\Delta = f_* - \min_{\substack{g(x) \leq \varepsilon_g \\ x \in Q}} f(x) = \min_{\substack{g(x) \leq 0 \\ x \in Q}} f(x) - \min_{\substack{g(x) \leq \varepsilon_g \\ x \in Q}} f(x),$$

то

$$P\left(f(\bar{x}^N) - f_* + \varepsilon_\Delta \geq 2(\varepsilon_f + \varepsilon_\Delta)\right) \leq \frac{E\left[f(\bar{x}^N)\right] - f_* + \varepsilon_\Delta}{2(\varepsilon_f + \varepsilon_\Delta)} \leq \frac{1}{2}.$$

Можно параллельно (независимо) запустить $\log_2(\sigma^{-1})$ траекторий метода. Обозначим $\bar{x}_{\min}^N$ тот из $\bar{x}^N$ на этих траекториях, который доставляет минимальное значение $f(\bar{x}^N)$. Из выписанного неравенства Маркова получаем, что имеет место неравенство

$$P\left(f(\bar{x}_{\min}^N) - f_* \geq 2\varepsilon_f + \varepsilon_\Delta\right) \leq \sigma.$$

К сожалению, этот подход требует малости $\varepsilon_\Delta$, что, вообще говоря, нельзя гарантировать из условий задачи.

Немного более аккуратные рассуждения (без новых идей) позволяют развязать во всех приведенных выше в п. 3 рассуждениях $\varepsilon_f$ и $\varepsilon_g$, допуская, что они могут выбираться независимо друг от друга. Детали мы вынуждены здесь опустить.

Основные приложения описанного подхода, это задачи вида

$$f(x) \to \min_{\max_{k=1,\ldots,m} \sigma_k(A_k^T x) \leq 0},$$

с разреженной матрицей

$$A = [A_1, ..., A_m]^T.$$

В частности, задачи вида

$$f(x) \to \min_{Ax \le b}$$

и приводящиеся к такому виду задачи

$$f(x) \to \min_{\substack{Ax \le b \\ Cx = d}}.$$

В этих задачах, как правило, выбирают $\|\cdot\| = \|\cdot\|_2$, $d(x) = \|x\|_2^2/2$. Подобно [10] можно попутно восстанавливать (без особых дополнительных затрат) и двойственные множители к этим ограничениям. Причем эта процедура позволяет сохранить дешевизну итерации даже в разреженном случае.

### 4. Примеры решения разреженных задач с использованием рандомизированного метода зеркального спуска

Начнем с известного примера [12], демонстрирующего практически все основные способы рандомизации, которые сейчас активно используются в самых разных приложениях. Рассмотрим задачу поиска левого собственного вектора $x$, отвечающего собственному значению 1, стохастической по строкам матрицы $P = \|p_{ij}\|_{i,j=1}^{n,n}$ (такой вектор называют вектором Фробениуса–Перрона, а с учетом контекста PageRank вектором). Изложение рандомизации, связанной с ускоренными покомпонентными методами мы опускаем, поскольку она не завязана на МЗС. Тем не менее, приведем ссылки на работы, в которых такой подход к поиску PageRank описан: замечание 5 [4] и пример 4 [13] (см. также замечания 10, 11 [14]).

Перепишем задачу поиска вектора PageRank следующим образом [3]

$$f(x) = \frac{1}{2}\|Ax\|_2^2 \to \min_{x \in S_n(1)},$$

где $S_n(1)$ – единичный симплекс в $\mathbb{R}^n$, $A = P^T - I$, $I$ – единичная матрица. Далее будем использовать обозначения $A^{\langle k \rangle}$ – $k$-й столбец матрицы $A$, $A_k$ – транспонированная $k$-я строка (то есть $A_k$ – это вектор) матрицы $A$. Следуя [12], воспользуемся для решения этой задачи рандомизированным МЗС со стохастическим градиентом[2]

$$\nabla_x f_k(x, \xi^k) = (P - I)^{\langle j(\xi^k) \rangle} - (P - I)^{\langle \xi^k \rangle},$$

где

$$\xi^k = i \text{ с вероятностью } x_i, \ i = 1, ..., n;$$
$$j(\xi^k) = j \text{ с вероятностью } p_{\xi^k j}, \ j = 1, ..., n.$$

Несложно проверить выполнение условия 1 п. 2, если генерирование использующихся вспомогательных случайных величин осуществляется независимо. В виду симплексных

---

[2] Сначала согласно вектору $x \in S_n(1)$ случайно разыгрываем один из столбцов матрицы $A = P^T - I$ (пусть это будет $\xi^k$-й столбец). Пользуясь тем, что столбцы матрицы $P^T$ сами представляют собой распределения вероятностей, независимо разыгрываем (согласно выбранному столбцу $\xi^k$) случайную величину, и выбираем соответствующий столбец матрицы $A^T = P - I$, из которого вычитаем $A^{T\langle \xi^k \rangle}$ (из-за наличия матрицы $I$), таким образом, конструируется несмещенная оценка градиента $\nabla f(x) = A^T A x$.

ограничений, естественно следовать при выборе прокс-структуры примеру 3 п. 2. Таким образом, можно оценить

$$\tilde{M}^2 = \max_{x \in S_n(1), \xi^k} \left\| \nabla_x f_k(x, \xi^k) \right\|_\infty^2 \le 4.$$

Даже в случае, когда матрица $P$ полностью заполнена амортизационная (средняя) стоимость одной итерации будет $\mathrm{O}(n)$ (вместо $\mathrm{O}(n^2)$, в случае честного расчета градиента). Таким образом, общее число арифметических операций будет $\mathrm{O}(n \ln n / \varepsilon^2)$.

К худшей оценке приводит другой способ рандомизации (рандомизации суммы [4]). Чтобы его продемонстрировать, перепишем исходную задачу следующим образом

$$f(x) = n \sum_{k=1}^n \frac{1}{n} \frac{1}{2} (A_k^T x)^2 \to \min_{x \in S_n(1)}.$$

Из такого представления следует, что можно определить стохастический градиент следующим образом

$$\nabla_x f_k(x, \xi^k) = n A_{\xi^k} A_{\xi^k j(x)},$$

где

$$\xi^k = i \text{ с вероятностью } 1/n, \; i = 1, ..., n;$$
$$j(x) = j \text{ с вероятностью } x_j, \; j = 1, ..., n.$$

Амортизационная (средняя) стоимость одной итерации будет по-прежнему $\mathrm{O}(n)$, но вот оценка $\tilde{M}^2$ получается похуже. Здесь мы имеем пример, когда $M^2$ и $\tilde{M}^2$ существенно отличаются – в действительности, можно вводить промежуточные условия, не такие жесткие, как условие 3, и получать более оптимистичные оценки вероятностей больших уклонений [4].

К сожалению, эти методы не позволяют полноценно воспользоваться разреженностью матрицы $P$, которая, как правило, имеет место. Собственно, этот пункт отчасти и будет посвящен тому, как можно сочетать рандомизацию и разреженность. В частности, если переписать задачу PageRank следующим образом

$$\|Ax\|_\infty \to \min_{x \in S_n(1)},$$

что равносильно (факт из теории неотрицательных матриц [15])

$$\max_{k=1,...,n} A_k^T x \to \min_{x \in S_n(1)},$$

то исходя из примера 3 (в детерминированном случае), можно получить следующую оценку [16] на общее число арифметических операций $\mathrm{O}(n \ln n / \varepsilon^2)$, при условии, что число элементов в каждой строке и столбце матрицы $P$ не больше $\mathrm{O}(\sqrt{n / \ln n})$. Здесь не использовалась рандомизация, а использовалась только разреженность матрицы $P$ (следовательно и $A$). По-сути, способ получения этой оценки всецело базируется на возможности организации эффективного пересчета субградиента функционала, подобно [1 – 3]. Далее мы распространим этот пример на более общий класс задач, и постараемся привнести в подход рандомизацию.

Итак, рассмотрим сначала класс задач с $Q$ из примера 1 или 2 п. 2

$$\max_{k=1,...,m} \sigma_k(A_k^T x) \to \min_{x \in Q}, \tag{7}$$

где $\sigma_k(\;)$ – выпуклые функции с константой Липшица равномерно ограниченной известным числом $\mathrm{M}$, (суб-)градиент каждой такой функции (скалярного аргумента) можно рассчитать за $\mathrm{O}(1)$. Введем матрицу

$$A = [A_1, ..., A_m]^T$$

и будем считать, что в каждом столбце матрицы $A$ не больше $s_m \le m$ ненулевых элементов, а в каждой строке – не больше $s_n \le n$. Заметим, что некоторую обременительность этим условиям создает требование, что в "каждом" столбце/строке. Это требования можно ослаблять, приближаясь к некоторым средним показателям разреженности (численные эксперименты в этой связи также проводились [3]), однако в данной работе для большей наглядности и строгости рассуждений мы ограничимся случаем, когда именно в каждом столбце/строке имеет место такая (или еще большая) разреженность.

Из работ [1 – 3] следует, что МЗС из примеров 1, 2 (в детерминированном случае) для задачи (7) будет требовать

$$O\left(\frac{\mathrm{M}^2 \max_{k=1,\ldots,m} \|A_k\|_2^2 R_2^2}{\varepsilon^2}\right)$$

итераций, где $R_2^2$ – квадрат евклидова расстояния от точки старта до решения, а одна итерация (кроме первой) будет стоить

$$O(s_n s_m \log_2 m).$$

И все это требует препроцессинг (предварительных вычислений, связанных с "правильным" приготовлением памяти) объема $O(m+n)$. Таким образом, в интересных для нас случаях общее число арифметических операций для МЗС из примеров 1, 2 будет

$$O\left(s_n s_m \log_2 m \frac{\mathrm{M}^2 \max_{k=1,\ldots,m} \|A_k\|_2^2 R_2^2}{\varepsilon^2}\right). \tag{8}$$

Постараемся ввести рандомизацию в описанный подход. Для этого осуществим дополнительный препроцессинг, заключающийся в приготовлении из векторов $A_k$ вектора распределения вероятностей. Представим

$$A_k = A_k^+ - A_k^-,$$

где каждый из векторов $A_k^+$, $A_k^-$ имеет не отрицательные компоненты. Согласно этим векторам приготовим память таким образом, чтобы генерирование случайных величин из распределений $A_k^+ / \|A_k^+\|_1$ и $A_k^- / \|A_k^-\|_1$ занимало бы $O(\log_2 n)$. Это всегда можно сделать [3]. Однако это требует хранение в "быстрой памяти" довольно большого количества соответствующих "деревьев". Весь этот препроцессинг и затраченная память будут пропорциональны числу ненулевых элементов матрицы $A$, что в случае huge-scale задач сложно осуществить из-за ресурсных ограничений. Тем не менее, далее мы будем считать, что такой препроцессинг можно осуществить, и (самое главное) такую память можно получить. Введем стохастический (суб-)градиент

$$\nabla_x f_k(x, \xi^k) = \|A_{k(x)}^+\|_1 e_{i(\xi^k)} - \|A_{k(x)}^-\|_1 e_{j(\xi^k)},$$

где

$$k(x) \in \operatorname*{Arg\,max}_{k=1,\ldots,m} \sigma_k(A_k^T x),$$

причем не важно, какой именно представитель Arg max выбирается;

$$e_i = (\underbrace{0,\ldots,0,1,0,\ldots,0}_{i});$$

$i(\xi^k) = i$ с вероятностью $A_{k(x)i}^+ / \|A_{k(x)}^+\|_1$, $i = 1,\ldots,n$;

$$j(\xi^k) = j \text{ с вероятностью } A^-_{k(x)j} \big/ \|A^-_{k(x)}\|_1, \ j=1,...,n.$$

Легко проверить выполнение условия 1 п. 2 (заметим, что $\nabla f(x) = A_{k(x)}$). Также легко оценить

$$\tilde{M}^2 \le M^2 \max_{k=1,...,m} \|A_k\|_1^2.$$

И получить из примеров 1, 2 следующую оценку числа итераций (ограничимся для большей наглядности сходимостью по математическому ожиданию, т.е. без оценок вероятностей больших уклонений)

$$O\left( \frac{M^2 \max_{k=1,...,m} \|A_k\|_1^2 R_2^2}{\varepsilon^2} \right)$$

Основная трудоемкость тут в вычислении $k(x)$. Однако, за исключением самой первой итерации можно эффективно организовать перерешивание этой задачи. Действительно, предположим, что уже посчитано $k(x^l)$, а мы хотим посчитать $k(x^{l+1})$. Поскольку согласно примерам 1, 2 $x^{l+1}$ может отличаться $x^l$ только в двух компонентах, то пересчитать $\max_{k=1,...,m} \sigma_k(A_k^T x^{l+1})$, исходя из известного $\max_{k=1,...,m} \sigma_k(A_k^T x^l)$, можно за (см., например, [1, 3]) $O(2s_m \log_2 m)$. Таким образом, общее ожидаемое число арифметических операций нового рандомизированного варианта МЗС из примеров 1, 2 для задачи (7) будет

$$O\left( s_m \log_2 m \frac{M^2 \max_{k=1,...,m} \|A_k\|_1^2 R_2^2}{\varepsilon^2} \right). \tag{9}$$

Для матриц $A$, все отличные от нуля элементы которых одного порядка, скажем $O(1)$, имеем

$$\max_{k=1,...,m} \|A_k\|_2^2 = s_n, \ \max_{k=1,...,m} \|A_k\|_1^2 = s_n^2.$$

В таком случае не стоит ожидать выгоды (формулы (8) и (9) будут выглядеть одинаково). Но если это условие (ненулевые элементы $A$ одного порядка) выполняется не очень точно, то можно рассчитывать на некоторую выгоду.

Рассмотрим теперь более общий класс задач, возникающих, например, при поиске равновесий в транспортных сетях [9]

$$\frac{1}{r} \sum_{k=1}^r \max_{l=a_k+1,...,b_k} \sigma_l(A_l^T x) \to \min_{x \in Q}, \tag{10}$$
$$0 = a_1 < b_1 = a_2 < b_2 = a_3 < ... < b_{r-1} = a_r < b_r = m.$$

Матрица $A$ и числа $s_n$, $s_m$ определяются аналогично. Привнося (при расчете стохастического градиента) к описанным выше двум подходам для задачи (7) сначала равновероятный (и независимый от других рандомизаций) выбор одного из слагаемых в этой сумме, получим соответствующие обобщения (для задачи (10)) оценок (8), (9), которые будут иметь точно такой же вид. Только матрица $A$ собирается теперь из всех слагаемых суммы (10).

Возвращаясь к примеру 3, заметим, что все описанные выше конструкции (в том числе, связанные с задачей (10)) можно перенести на этот пример, в случае, когда

$$\sigma_k(A_k^T x) = A_k^T x - b_k.$$

При этом

$$R_2^2 \to \ln n, \ M = 1 \ (\text{для } (8), (9))$$

$$\max_{k=1,\ldots,m} \|A_k\|_2^2 \to \max_{\substack{i=1,\ldots,m \\ j=1,\ldots,n}} |A_{ij}|, \ \min\{s_n s_m \log_2 m, n\} \to \max\{s_n s_m \log_2 m, n\} \ (\text{для (8)})$$

$$\max_{k=1,\ldots,m} \|A_k\|_1^2 \to \max_{k=1,\ldots,m} \|A_k\|_1^2, \ s_m \log_2 m \to \max\{s_m \log_2 m, n\} \ (\text{для (9)}).$$

Собственно, примера PageRank, изложенный в начале этого пункта, как раз подходил под применение оценки (8).

С помощью п. 3 все написанное выше переносится и на задачи вида

$$f(x) \to \min_{\substack{Ax \le b \\ Cx = d}},$$

с разреженными матрицами. Такие задачи играют важную роль, например, при проектировании механических конструкций [2] (Truss topology design). Мы не будем здесь приводить соответствующие рассуждения, поскольку они достаточно очевидны, и заинтересованный читатель сможет осуществить отмеченное обобщение самостоятельно.

## Литература